\documentclass[12pt]{amsart}
\usepackage{amsmath,amscd,amssymb,amsfonts}
\newcommand{\q}{\quad}
\newcommand{\h}{\hbox}
\newcommand{\nin}{\noindent}
\newcommand{\bl}{\bigl}
\newcommand{\br}{\bigr}
\newcommand{\bs}{\par\bigskip}
\newcommand{\ms}{\par\medskip}
\newcommand{\sk}{\par\smallskip}
\newcommand{\md}{\C\{\!\{\dti\}\!\}}
\newcommand{\ssb}{\raise.2ex\h{${\scriptscriptstyle\bullet}$}}
\newcommand{\ssc}{\,\raise.2ex\h{${\scriptstyle\circ}$}\,}
\newcommand{\sotim}{\h{$\otimes$}}

\newcommand{\mtim}{\h{$\times$}}
\newcommand{\mopl}{\h{$\bigoplus$}}
\newcommand{\mcap}{\h{$\bigcap$}}
\newcommand{\mcup}{\h{$\bigcup$}}
\newcommand{\msum}{\h{$\sum$}}
\newcommand{\mprod}{\h{$\prod$}}
\newcommand{\lan}{\langle}
\newcommand{\ran}{\rangle}
\newcommand{\ges}{\geqslant}
\newcommand{\les}{\leqslant}
\newcommand{\al}{\alpha}
\newcommand{\A}{{\mathcal A}}
\newcommand{\B}{{\mathcal B}}
\newcommand{\C}{{\mathbf C}}
\newcommand{\dd}{\partial}
\newcommand{\dt}{\partial_t}
\newcommand{\dti}{\partial_t^{-1}}
\newcommand{\D}{{\mathcal D}}

\newcommand{\De}{\Delta}
\newcommand{\E}{{\mathcal E}}
\newcommand{\EE}{\widehat{\mathcal E}}
\newcommand{\Eh}{\widehat{E}}
\newcommand{\ep}{\varepsilon}
\newcommand{\fb}{{}\,\overline{\!f}{}}
\newcommand{\F}{{\mathcal F}}
\newcommand{\G}{{\mathcal G}}
\newcommand{\GG}{\widetilde{\mathcal G}}
\newcommand{\Hc}{{\mathcal H}}
\newcommand{\Htc}{\widetilde{\mathcal H}}
\newcommand{\Hhc}{\widehat{\mathcal H}}
\newcommand{\I}{{\mathcal I}}
\newcommand{\J}{{\mathcal J}}
\newcommand{\K}{{\mathcal K}}
\newcommand{\Kb}{K^{\ssb}}

\newcommand{\Kh}{\widehat{K}}
\newcommand{\Kch}{\widehat{\mathcal K}}
\newcommand{\La}{\Lambda}
\newcommand{\la}{\lambda}
\newcommand{\Lc}{{\mathcal L}}
\newcommand{\Lb}{{}\,\overline{\!{\mathcal L}}{}}
\newcommand{\LL}{{\mathcal L}^{\ssb}}
\newcommand{\M}{{\mathcal M}}
\newcommand{\MM}{{}\,\overline{\!{\mathcal M}}{}}

\newcommand{\NN}{{\mathbf N}}
\newcommand{\PP}{{\mathbf P}}
\newcommand{\pib}{\overline{\pi}}
\newcommand{\OO}{{\mathcal O}}
\newcommand{\OS}{{\mathcal O}_S}
\newcommand{\OX}{{\mathcal O}_X}
\newcommand{\OXb}{\OO_{\X}}
\newcommand{\Q}{{\mathbf Q}}
\newcommand{\R}{{\mathbf R}}
\newcommand{\Rh}{\widehat{R}}

\newcommand{\V}{{\mathcal V}}
\newcommand{\X}{{}\,\overline{\!X}{}}
\newcommand{\XX}{{\mathcal X}}
\newcommand{\xit}{\widetilde{\xi}}
\newcommand{\Z}{{\mathbf Z}}
\newcommand{\DR}{{\rm DR}}
\newcommand{\an}{{\rm an}}
\newcommand{\Diff}{{\rm Diff}}
\newcommand{\Ker}{{\rm Ker}}

\newcommand{\Hom}{{\mathcal H}om}

\newcommand{\Gr}{\h{\rm Gr}}
\newcommand{\1}{{\hskip1pt}}
\newcommand{\simto}{\buildrel\sim\over\longrightarrow}

\newcommand{\into}{\hookrightarrow}
\newcommand{\ilim}{\rlap{\raise-5.5pt\h{$\,\rightarrow$}}{\rm lim}}
\newcommand{\plim}{\rlap{\raise-5.5pt\h{$\,\leftarrow$}}{\rm lim}}
\begin{document}
\title[Kontsevich's conjecture]
{Kontsevich's conjecture on an algebraic formula\\
for vanishing cycles of local systems}
\author[C. Sabbah]{Claude Sabbah}
\address{UMR 7640 du CNRS, Centre de Math\'ematiques Laurent Schwartz,
Ecole polytechnique, F-91128 Palaiseau cedex, France}
\author[M. Saito]{Morihiko Saito}
\address{RIMS Kyoto University, Kyoto 606-8502 Japan}
\thanks{
The first named author is supported by the grant ANR-08-BLAN-0317-01
of the Agence nationale de la recherche.
The second named author is partially supported by Kakenhi 21540037.}
\begin{abstract}
For a local system and a function on a smooth complex algebraic
variety, we give a proof of a conjecture of M. Kontsevich on a
formula for the vanishing cycles using the twisted de Rham complex of
the formal microlocalization of the corresponding locally free sheaf
with integrable connection having regular singularity at infinity.
We also prove its local version, which may be viewed as a natural
generalization of a result of E. Brieskorn in the isolated singularity
case. We then generalize these to the case of the de Rham complexes of
regular holonomic D-modules where we have to use the tensor product
with a certain sheaf of formal microlocal differential operators
instead of the formal completion.
\end{abstract}
\maketitle
\bs\bs
\centerline{\bf Introduction}
\bs\nin
Let $X$ be a smooth complex algebraic variety, and
$f\in\Gamma(X,\OX)$.
Let $L$ be a $\C$-local system on $X$.
For each $c\in\C$, we have the vanishing cycle sheaf complex
$\varphi_{f-c}\,L$ on $X_c:=f^{-1}(c)$ with the monodromy $T_c$, as is
defined by Deligne [De3] (see also Notation below).
Let $S=\C$ as a complex algebraic variety endowed with the natural
coordinate $t$. Then $f$ is identified with a morphism $f:X\to S$,
and $f^*t$ with the function $f$. Set $\dt:=\dd/\dd t$, and define
$$\Rh:=\C[[\dti]],\q\Kh:=\C((\dti))=\Rh[\dt].$$
Let $\Kh\lan t\ran=\mopl_{i\in\NN}\,\Kh t^i$.
This has a structure of a noncommutative ring by $\dt t-t\dt=1$.
For a finite dimensional $\C$-vector space $V$ with an automorphism
$T$ and for $c\in\C$, define a $\Kh\lan t\ran$-module by
$$\EE(V,T)_c:=V((\dti))=V[[\dti]][\dt],$$
with action of $t$ given by
$$t(v\,\dt^j)=\Bigl(-\frac{\log T}{2\pi i}-j\Bigr)v\,\dt^{j-1}+
c\,v\,\dt^j\q\h{for}\,\,\,v\in V,\,j\in\Z,$$
where the real part of the eigenvalues of $(2\pi i)^{-1}\log T$ are
contained in $[0,1)$.
\sk
Corresponding to the local system $L$, we have a locally free
$\OX$-module $\M$ with an integrable connection $\nabla$ having
regular singularities at infinity [De1].
Set $u:=\dti$.
We have a twisted de Rham complex
$$\M((u))^f\sotim_{\OX}\Omega_X^{\ssb}:=\DR_X\bl(\M((u))^f\br)
[-\dim X],$$
where $\M((u))^f:=\M[[u]][u^{-1}]$ and the differential of
$\DR_X\bl(\M((u))^f\br)$ is given by
$$\nabla-u^{-1}\sotim df.$$
Note that $\M((u))^f$ has a natural action of $\Kh=\C((u))$
together with an action of $t$ defined by
$$t(m\,\dt^j)=fm\,\dt^j-jm\,\dt^{j-1}
\q\h{for}\,\,\,m\in\M,\,j\in\Z.$$

In this paper we give a proof of the following which was conjectured
by M.~Kontsevich as a possible answer to a question raised in [KoSo]
concerning a definition of vanishing cycles for functions on smooth formal
schemes. We also simplify some arguments of a different approach in [Sab2].

\ms\nin
{\bf Theorem~1.} {\it For $k\in\Z$, there are canonical
isomorphisms as $\Kh\lan t\ran$-modules
$$H^k\bl(X,\M((u))^f\sotim_{\OX}\Omega_X^{\ssb}\br)=
\mopl_{c\in\C}\,\EE\bl(H^{k-1}(X_c,\varphi_{f-c}\,L),T_c^{(k-1)}\br)_c,$$
where $T_c^{(k-1)}$ is induced by $T_c$.}
\sk
Note that the right-hand side is a finite direct sum since
$\varphi_{f-c}\,L=0$ except for a finite number of $c$.
Viewed as a formal meromorphic connection for the variable $u$, the above
decomposition is a special case of the classical Levelt-Turrittin theorem.
Let $\V^k$ denote the left-hand side of the formula. Theorem~1 means that, in order to calculate $(H^{k-1}(X_c,\varphi_{f-c}\,L),T_c^{(k-1)})$,
it is enough to find a finite dimensional $\C$-vector subspace $V^k$ of
$\V^k$ such that
$$\Kh\sotim_{\C}V^k\simto\V^k,\q t\,V^k\subset\Rh\,V^k,$$
and moreover $A_0$ is semisimple and any two eigenvalues of $A_1(c,c)$ do
not differ by an integer for each $c$.
Here $\sum_{i\ges 0}A_i\,\dt^{-i}$ is the expansion of the action of $t$
with $A_i\in{\rm End}_{\C}(V^k)$, $V^k=\mopl_c\,V^k_c$ is the eigenspace
decomposition by the action of $A_0$, and $A_1(c,c)$ is the
$(V^k_c,V^k_c)$-component of $A_1$. In this case, we have the following
isomorphism (see Remark~(3.6) below)
$$(H^{k-1}(X_c,\varphi_{f-c}\,L),T_c^{(k-1)})=(V_c^k,\exp(-2\pi iA_1(c,c))).$$
In case $f$ has only isolated singular points, this calculation is
essentially equivalent to the one given by E.~Brieskorn [Br] (see (3.5)
below). We have a Cech calculation in the general case (see (3.4)
below). In special cases, we have the following.

\ms\nin
{\bf Corollary~1.} {\it If $X$ is affine, then the formula in Theorem~$1$
holds with the left-hand side replaced by the $k$-th cohomology of the
complex whose $q$-th component is
$$\Gamma\bl(X,\M\sotim_{\OO_X}\Omega_X^q\br)((u))^f,$$
where the differential is induced by $\nabla-u^{-1}df\wedge$.}

\ms\nin
{\bf Corollary~2.} {\it In case $L=\C_X$ and $X$ is an affine open
subvariety of $\C^n$, let $A_X$ be the affine ring of $X$ with
$x_1,\dots,x_n$ the coordinates of $\C^n$.
Set $\dd_{x_i}:=\dd/\dd x_i$, and $f_i:=\dd f/\dd x_i$. The formula
in Theorem~$1$ holds with the left-hand side replaced by
$$H^kK^{\ssb}\bl(A_X((u));\dd_{x_i}-u^{-1}f_i\,(i\in[1,n]\br),$$
which is the cohomology of the Koszul complex for the action of
$\dd_{x_i}-u^{-1}f_i\,(i\in[1,n])$ on $A_X((u))$.}

\ms
What is interesting in Theorem~1 and its corollaries is that we have {\it only}
the contributions of the vanishing cycles of $f$ inside $X$ and {\it no
contributions of the vanishing cycles at infinity} on the right-hand
side of the equality. This might be rather surprising, since the
cohomological direct images $\Hc^jf_*\M$ of $\M$ as $\D$-modules
are defined by
$$\Hc^jf_*\M:=\Hc^jf_{\ssb}(\M[\dt]^f\sotim_{\OX}\Omega_X^{\ssb}
[\dim X]),$$
where $f_{\ssb}$ denotes the sheaf-theoretic direct image, and we have
$\dt=u^{-1}$ so that $\M((u))^f$ is isomorphic to the {\it formal
microlocalization} of $\M[\dt]^f$. Note that the vanishing cycles of the
cohomological direct images $\Hc^jf_*\M$ have the contributions from
the vanishing cycles at infinity in general. So these imply that one
cannot apply the formal microlocalization after taking the direct image
of $\M$ by $f$ as a $\D$-module.
Note also that the vanishing cycle functor does not necessarily commute
with the direct image under a nonproper morphism $f$ because of the
problem of singularities at infinity of $f$.
\sk
Let $j:X\into\X$ be a smooth compactification such that $f$ induces
a morphism $\fb:\X\to\PP^1$. Set $\MM=j_*\M$.
Using GAGA we can show the following (see also [Sab2]).

\ms\nin
{\bf Proposition~1.} {\it For $k\in\Z$, there are canonical
isomorphisms of $\Kh\lan t\ran$-modules}
$$H^k\bl(\X,\MM((u))^f\sotim_{\OO_{\X}}\Omega_{\X}^{\ssb}\br)\simto
H^k\bl(\X^{\an},\MM^{\an}((u))^f\sotim_{\OO_{\X^{\an}}}
\Omega_{\X^{\an}}^{\ssb}\br).$$
\ms
Note that the noncommutativity of cohomology and inductive limit
does not cause a problem for the construction of the above canonical
morphism by using the continuous morphism $\rho:\X^{\an}\to\X$,
see (3.1).
By Proposition~1, the proof of Theorem~1 can be reduced to local
analytic calculations on $\X^{\an}$.
So we will consider only the underlying complex manifolds,
and the upperscript $^{\an}$ will be omitted from now on in this paper
(except for (2.1), (2.3) and (3.1)).
Moreover it is enough to work mainly on $X$ by Proposition~4 below.
\sk
We also have the twisted de Rham complex
$\M[\dt,\dti]^f\sotim_{\OX}\Omega_X^{\ssb}$ associated with the
algebraic localization $\M[\dt,\dti]^f$.
Consider its cohomology sheaves
$$\Hc^j\bl(\M[\dt,\dti]^f\sotim_{\OX}\Omega_X^{\ssb}\br).$$
These are constructible sheaves of Gauss-Manin systems in one-variable
(up to the division by $\C[t]$ for $j=1$), see [BaSa].
We show that the topology on their stalk at each $x\in X$ induced
by the $u$-adic topology on $\M[u]^f$ coincides with the one induced
by the $\dti$-adic topology on the Brieskorn lattices of the
Gauss-Manin systems by using the finiteness theorem on the order of
the torsion of Brieskorn modules proved in loc.~cit.
We then get the following.

\ms\nin
{\bf Proposition~2.} {\it There is a canonical quasi-isomorphism of
complexes of $f^{-1}\EE_S$-modules
$$f^{-1}\EE_S\sotim_{f^{-1}\D_S}\bl(\M[\dt]^f\sotim_{\OX}
\Omega_X^{\ssb}\br)\simto
\M((\dti))^f\sotim_{\OX}\Omega_X^{\ssb},$$
where $\M[\dt]^f$ is the direct image of $\M$ as an analytic
$\D$-module by the graph embedding of $f$, and $\EE_S$ is the
ring of formal microdifferential operators on the complex
manifold $S=\C$, i.e. $\EE_S=\OS((\dti))$ forgetting the
multiplicative structure.}
\ms
Here we also use the assertion that the filtration $V$ of Kashiwara
and Malgrange on $\M[\dt]^f$ induces the $V$-filtration on each
stalk of the constructible cohomology sheaves of the Gauss-Manin
systems by taking the relative de Rham complex, see (1.3.3) below.
\sk
Using Proposition~2 we get the following local analytic version of
Theorem~1 for local systems $L$ and holomorphic functions $f$ on
complex manifolds $X$ where $\M=\OX\sotim_{\C}L$ in the analytic
case.

\ms\nin
{\bf Theorem~2.} {\it For $k\in\Z$, there are canonical isomorphisms
of analytic constructible sheaves of $\Kh\lan t\ran$-modules on $X_c$
$$\Hc^{k+1}\bl(\M((u))^f\sotim_{\OX}\Omega_X^{\ssb}\br)\big|_{X_c}=
\EE\bl(\Hc^k\varphi_{f-c}\,L,T_c^{(k-1)}\br)_c,$$
where $\EE(*)_c$ on the right-hand side is naturally extended to the
case of constructible sheaves.}
\bs
This would not be very surprising to the specialists in view of
[BaSa] and also [Sai2] (see Remark~(1.7)(i) and Remark~(3.5) below).
If $f$ has only isolated singular points, then it is essentially
equivalent to the calculation of Brieskorn [Br] explained after
Theorem~1.

\ms
For the proof of Theorem~1, we have a passage from local to global,
and need the following.

\ms\nin
{\bf Proposition~3.} {\it The filtration $V$ of Kashiwara and
Malgrange on $\M[\dt]^f$ along $t=c$ induces the filtration $V$ on
the global analytic cohomology $H^k\bl(X_c,\M[\dt]^f\sotim_{\OX}
\Omega_X^{\ssb}|_{X_c}\br)$ in a strict way, i.e.\ $\Gr_V^{\al}$
commutes with the global cohomology.}
\ms
This is proved by using a filtered spectral sequence which
is a priori defined in an abelian category containing the exact
category of filtered vector spaces, see (3.2) below.
\sk
For the proof of Theorem~1, there still remains the following
key proposition.

\ms\nin
{\bf Proposition~4.} {\it We have a canonical analytic
quasi-isomorphism}
$$\MM((u))^f\sotim_{\OXb}\Omega_{\X}^{\ssb}\simto
\R j_*\bl(\M((u))^f\sotim_{\OX}\Omega_X^{\ssb}\br).$$
\sk
Note that the right-hand side is supported on a union of finite
number of $\X_c\,\,(c\in\C)$ by Theorem~2, and this implies the
vanishing of the restriction of the left-hand side to $\X_{\infty}$.
The proof of Proposition~4 is reduced to the case where the union of
$D:=\X\setminus X$ and the singular $X_c$ is a divisor with normal
crossings on $\X$.
Indeed, for a proper morphism of complex manifolds $\pi:X'\to X$,
the canonical morphism $L\to\R\pi_*\pi^*L$ splits by applying this
morphism to the dual of $L$ and using Verdier duality.
Here we use the direct image of differential complexes rather
than that of $\D$-modules to simplify the argument.
For the proof of Proposition~4 in the normal crossing case,
we prove a sufficient condition for the commutativity of the global
section functor and the inductive limit in (1.8), and apply this to
the filtration $G$ used in the proof of Theorem~2.
\sk
In this paper we also show the following.

\ms\nin
{\bf Theorem~3.} {\it Theorems~$1$ and $2$ remain valid by replacing
respectively $\M$, $L$, and $\M((u))^f$ with a bounded complex of algebraic
regular holonomic $\D_X$-modules $\M^{\ssb}$, $\DR_X(\M^{\ssb})[-\dim X]$,
and $\M^{\ssb}\sotim_{\OX}\OX((u))^f$.}
\ms
Here the replacement of $\M((u))^f$ by $\M^{\ssb}\sotim_{\OX}\OX((u))^f$
is quite essential. Indeed, we cannot get a correct result if we apply
the same definition as in the local system case.
This is a quite subtle point of the theory, and is related to the
noncommutativity of inductive limit and projective limit.
Using canonical isomorphisms, the proof of Theorem~3 can be reduced
to the case $\M$ is a locally free $\OX(*E)$-module having a regular
singular connection. Here $E$ is a divisor on $X$, called the
{\it interior} divisor.
Note also that, setting $\MM=j_*\M$ with $j:X\into\X$ as above, the
construction of $\MM((u))^f$ along the interior divisor $E$ is quite
different from the one along the exterior divisor $D:=\X\setminus X$,
see Example~(2.2).
We prove Theorem~3 after showing Theorems~1 and 2 in the local system
case since many readers would be interested mainly in the original
conjecture of Kontsevich in this special case.
\sk
We thank the referee for useful comments to improve the paper.
\sk
In Section~1 we review some basics about the sheaves of Gauss-Manin
systems and Brieskorn modules, and then prove Theorem~2 after showing
Proposition~2.
In Section~2 we recall some basics of formal microlocalization, and
prove Proposition~4 by reducing to the normal crossing case.
In Section~3 we prove Theorem~1 after showing Propositions~1 and 3.
In Section~4 we show how to generalize the arguments in the previous
sections in order to prove Theorem~3.

\ms\nin
{\bf Notation.}
The nearby and vanishing cycle functors in [De3] are denoted
respectively by $\psi_f$ and $\varphi_f$.
In case $\psi_f\F$ is a (shifted) perverse sheaf and $\F$ has
$\C$-coefficients, $\psi_{f,\la}\F$ denotes the $\la$-eigen-subsheaf
of $\psi_f\F$ for the semisimple part $T_s$ of the monodromy $T$,
i.e. $\psi_{f,\la}\F:=\Ker(T_s-\la)\subset\psi_f\F$, and similarly
for $\varphi_{f,\la}\F$.

\bs\bs
\centerline{\bf 1. Sheaves of Gauss-Manin systems and Brieskorn
modules}
\bs\nin
In this section we review some basics about the sheaves of Gauss-Manin
systems and Brieskorn modules, and then prove Theorem~2 after showing
Proposition~2.
Here we treat only complex manifolds and analytic sheaves.

\ms\nin
{\bf 1.1.~Regular holonomic $\D_{S,0}$-modules.}
Let $S=\C$ as a complex manifold. Let $t$ be the coordinate.
Set $\dt=\dd/\dd t$.
For a regular holonomic left $\D_{S,0}$-module $M$ and $\al\in\C$,
define
$$M^{\al}:=\Ker\bl((t\dt-\al)^k:M\to M\br)\q(k\gg 0).$$
The actions of $t$ and $\dt$ on $M$ induce the morphisms for
any $\al\in\C$
$$t:M^{\al}\to M^{\al+1},\q \dt:M^{\al+1}\to M^{\al},
\leqno(1.1.1)$$
which are bijections for $\al\ne-1$.
We have the canonical inclusions
$$\mopl_{\al\in\C}\,M^{\al}\subset M\subset\mprod_{\al\in\C}\,M^{\al}.
\leqno(1.1.2)$$
The last inclusion follows from the fact that
$\bigoplus_{\al\in\C}M^{\al}$ generates $M$ over $\OO_{S,0}=\C\{t\}$.
The latter can be proved by using the classical theory of
differential equations of one variable with regular singularities,
see e.g.\ [De1].
(Note that (1.1.2) does not hold in the higher dimensional case
unless $M$ is isomorphic to its Verdier specialization along the
hypersurface.)
\sk
Take an extension of $M$ to a coherent $\D_S$-module defined on a
sufficiently small open neighborhood of $0\in S$, which is also
denoted by $M$. We then get a perverse sheaf $\DR_S(M)$, replacing $S$
with a neighborhood of $0$. By [Kas1], [Mal], there are canonical
isomorphisms of $\C$-vector spaces for $\la=\exp(-2\pi\al)$
$$M^{\al}=\begin{cases}\psi_{t,\la}\DR_S(M)[-1]&
\h{if}\,\,\,\al\notin\Z_{\les-1},\\
\varphi_{t,\la}\DR_S(M)[-1]&\h{if}\,\,\,\al\notin\Z_{\ges 0},\end{cases}
\leqno(1.1.3)$$
where $\Z_{\les k}:=\{i\in\Z\mid i\les k\}$, and similarly for
$\Z_{\ges k}$.
(For $\psi_{t,\la}$, $\varphi_{t,\la}$, see Notation at the end of the
introduction.)
Set
$$\La:=\{\al\in\C\mid {\rm Re}\,\al\in[0,1)\}.
\leqno(1.1.4)$$
We have canonical isomorphisms of $\C$-vector spaces
$$\psi_t\DR_S(M)[-1]=\mopl_{\al\in\La}\,M^{\al},\q
\varphi_t\DR_S(M)[-1]=\mopl_{\al\in\La-1}\,M^{\al}.
\leqno(1.1.5)$$
Indeed, let $V$ denote the $V$-filtration of Kashiwara and Malgrange
indexed decreasingly by $\Z$ so that the eigenvalues of the action
of $t\dt$ on $\Gr_V^iM$ are contained in $\La+i$,
see loc.~cit. Then
$$V^iM=M\cap\mprod_{{\rm Re}\,\al\ges i}\,M^{\al}.
\leqno(1.1.6)$$
This implies the canonical splittings (using the coordinate $t$)
for $i>j$
$$V^jM=\mprod_{j\les{\rm Re}\,\al<i}\,M^{\al}\oplus V^iM
\leqno(1.1.7)$$

\ms\nin
{\bf 1.2.~Microlocalizations.}
Let $\E_S$ be the ring of microdifferential operators on $P^*S\,(=S)$,
and $\EE_S$ be the ring of formal ones, see [SKK], [Kas2].
Here the
projective cotangent bundle $P^*S$ is identified with $S$ since
$\dim S=1$.
Forgetting the multiplicative structure, we have
$$\EE_S=\OS((\dti))\,(:=\OS[[\dti]][\dt]).
\leqno(1.2.1)$$
It is known that $\E_S$ and $\EE_S$ are flat over $\D_S$,
see loc.~cit.
So the tensor with $\E_{S,0}$ or $\EE_{S,0}$ over $\D_{S,0}$ is
an exact functor. Set
$$\aligned &R:=\md\,\bl(:=\bl\{\msum_{i\in\NN}\,a_i\dt^{-i}\,
\big|\,\msum_{i\in\NN}|a_i|r^i/i!<\infty\,\,\,\h{for some}\,\,r>0
\br\}\br),\\
&K:=R[\dt],\q\Rh=\C[[\dti]],\q\Kh:=\Rh[\dt]=\C((\dti)).\endaligned
\leqno(1.2.2)$$

For $M$ as in (1.1), we then get the canonical isomorphisms
$$\aligned \E_{S,0}\sotim_{\D_{S,0}}M&\simto
\mopl_{\al\in\La}K\sotim_{\C}M^{\al},\\
\EE_{S,0}\sotim_{\D_{S,0}}M&\simto
\mopl_{\al\in\La}\Kh\sotim_{\C}M^{\al}.\endaligned
\leqno(1.2.3)$$
To define the action of $\E_{S,0}$, $\EE_{S,0}$ on the right-hand
side, we use the expression
$$P=\msum_{j\in\Z}\,\msum_{i=0}^{m-j}\,a_{i,j}(t\dt)^i\dt^j\q
(a_{i,j}\in\C)\q\h{for}\,\,\,\,P\in F_m\E_{S,0}\,\,\,\h{or}
\,\,\,F_m\EE_{S,0},$$
where $F$ is the filtration by the order of $\dt$.
This induces the canonical morphisms in (1.2.3) by using (1.1.2).
These are isomorphisms by reducing to the case $M$ simple, i.e.
$M=\D_{S,0}/\D_{S,0}Q$ with $Q$ either $\dt$ or $t$ or
$t\dt-\al\,\,(\al\in\La\setminus\{0\})$.

\ms\nin
{\bf 1.3.~Sheaves of Gauss-Manin systems.}
Let $f$ be a nonconstant holomorphic function on a complex
manifold $X$.
Let $i_f:X\to X\times S$ be the graph embedding by $f$.
Let $\B_f$ be the direct image of $\OX$ by $i_f$ as a left
$\D_X$-module. Then
$$\B_f=\OX[\dt]^f,$$
where the sheaf-theoretic direct image by $i_f$ (i.e. the zero
extension) is omitted to simplify the notation, and $^f$ means that
the actions of $t$ and vector fields $\xi$ on $X$ are twisted so that
$$\aligned \xi(g\,\dt^i)&=\xi g\,\dt^i-(\xi f)\,g\,\dt^{i+1}
\q\h{for}\,\,g\in\OX,\\
t(g\,\dt^i)&=fg\,\dt^i-ig\,\dt^{i-1}.\endaligned$$
(Here the delta function $\delta(f-t)$ is actually omitted after
$\dt^i$.)
Set $X_0=f^{-1}(0)$, and define
$$\K_f^{\ssb}=\DR_{X\times S/S}(\B_f)[-\dim X]|_{X_0}.$$
This is a complex whose $i$-th term is $\Omega_X^i[\dt]|_{X_0}$
and whose differential is given by
$$d(\omega\,\dt^i)=(d\omega)\,\dt^i-(df\wedge\omega)\,\dt^{i+1}\q
\h{for}\,\,\omega\in\Omega_X^i.$$
Set
$$\G_f^i:=\Hc^i\K_f^{\ssb},$$
where $\Hc^i$ denotes the cohomology sheaf of a sheaf complex.
We call $\G_f^i$ the {\it sheaf of Gauss-Manin systems}
associated to $f$.
It is a sheaf of regular holonomic left $\D_{S,0}$-modules.
\sk
For $i=1$, we see that $\omega_0:=df\in\G_f^{1}$ is annihilated by
$\dt$ (i.e. $df\,\dt=0$ in $\G_f^{1})$, considering the image of
$1$ by $d$. So $\D_{S,0}\,\omega_0=\C\{t\}\,\omega_0$.
This is a free $\C\{t\}$-module of rank $1$ using inductively
$[\dt,t^i]=it^{i-1}$.
Define the {\it sheaf of reduced Gauss-Manin systems} $\GG_f^i$ by
$$\GG_f^i:=\begin{cases}\G_f^i&\h{if}\,\,\,i\ne 1,\\
\G_f^{1}/\C\{t\}\omega_0&\h{if}\,\,\,i=1.\end{cases}
\leqno(1.3.1)$$

It is well known that the stalks $\G_{f,x}^i$ are regular holonomic
$\D_{S,0}$-modules, and in the notation of (1.1), we have moreover
canonical isomorphisms of constructible sheaves for
$\la=\exp(-2\pi i\al)$
$$(\G_f^{i+1})^{\al}=\begin{cases}\Hc^i\psi_{f,\la}\C_X&
\h{if}\,\,\,\al\notin\Z_{\les 0},\\
\Hc^i\varphi_{f,\la}\C_X&\h{if}\,\,\,\al\notin\Z_{\ges 1}.\end{cases}
\leqno(1.3.2)$$
(For $\psi_{f,\la}\C_X$, $\varphi_{f,\la}\C_X$, see Notation at the
end of the introduction.)
\sk
Let $V$ denote the $V$-filtration of Kashiwara [Kas1] and
Malgrange [Mal] on $\B_f$ such that the roots of the minimal
polynomial of the action of $t\dt$ on $\Gr_V^i\B_f$ are contained
in $\La+i$ where $\La$ is as in (1.1.4).
Then (1.3.2) is reduced to the assertion that the $V$-filtration on
$\G_{f,x}^i$ is strictly induced by the $V$-filtration on $\B_f$,
i.e.\
$$V^{\al}\G_{f,x}^i={\rm Im}(\Hc^iV^{\al}\K_{f,x}^{\ssb}\into
\Hc^i\K_{f,x}^{\ssb}),
\leqno(1.3.3)$$
where ``strictly" means the injectivity of the morphism from
$\Hc^iV^{\al}\K_{f,x}^{\ssb}$.
The proof of (1.3.3) is further reduced to the assertion that
the induced filtration on $\G_{f,x}^i$ consists of finite
$\OO_{S,0}$-modules since the other conditions of the $V$-filtration
can be proved easily by using the action of $t\dt$ together with
this property.
Then the assertion follows from [Sai3], Prop.~3.4.8 and Lemma~3.4.9.
For the proof of (1.3.2) we also use the inclusion of Milnor tubes
for $x\in X_0$ and $x'\in X_0$ sufficiently near $x$ to show the
compatibility of the isomorphism (1.3.2) with the sheaf structure.
\sk
By definition there is a distinguished triangle
$$\C_{X_0}\to\psi_f\C_X\to\varphi_f\C_X\to,
\leqno(1.3.4)$$
In particular, $\Hc^i\psi_{f,\la}\C_X=\Hc^i\varphi_{f,\la}\C_X$
if $i\ne 0$ or $\la\ne 1$.
\sk
By (1.3.2) and (1.3.4), $\GG_f^i$ is the microlocalization of $\G_f^i$,
i.e.
$$\GG_f^i=\E_{S,0}\sotim_{\D_{S,0}}\G_f^i\q\h{for any}\,\,\,i\in\Z,
\leqno(1.3.5)$$
and the $\GG_f^i$ are finite dimensional $K$-vector spaces,
where $K$ is as in (1.2.2).

\ms\nin
{\bf 1.4.~Sheaves of Brieskorn modules.}
With the notation of (1.3), let $\A_f^{\ssb}$ be the complex defined
$$\A_f^i:=\Ker\bl(df\wedge :\Omega_X^i\to\Omega_X^{i+1}\br)\big|_{X_0},$$
with differential induced by $d$. There is a natural inclusion
$$\A_f^{\ssb}\to\K_f^{\ssb}.
\leqno(1.4.1)$$
Let
$$\Hc_f^i=\Hc^i\A_f^{\ssb}.$$
Note that we have to take the restriction to $X_0$ in order to
define the action of $\dti$ on $\Hc^{1}\A_f^{\ssb}$, see [BaSa].
We will call $\Hc_f^i$ the {\it sheaf of Brieskorn modules}
associated to $f$.
The inclusion (1.4.1) induces the canonical morphism
$$\Hc_f^i\to\G_f^i.
\leqno(1.4.2)$$

Since the differential is $f^{-1}\OS$-linear, there is a natural
structure of a $\C\{t\}$-module on $\Hc_f^i$.
We have the action of $\dti$ on $\Hc_f^i$ defined by
$$\dti\omega = df\wedge\eta\quad
\h{with}\quad d\eta =\omega.
\leqno(1.4.3)$$
This action is well defined.
For $i\ne 1$, the ambiguity of $\eta$ is given by $d\eta'$,
but we have $df\wedge d\eta'=-d(df\wedge\eta')$.
For $i=1$, the ambiguity of $\eta\in\OX$ is given by $\C$, and we
have a canonical choice of $\eta$ assuming that the restriction of
$\eta$ to $X_0$ vanishes (this is allowed because
$df\wedge d\eta = 0$).
\sk
For $i=1$, we see that $\omega_0:=df\sotim 1$ in (1.3)
belongs to $\Hc_f^{1}$.
Here we can easily verify that $\Hc_f^{1}\to\G_f^{1}$ is injective,
see the proof of Th.~2.2 in [BaSa].
We define the {\it sheaf of reduced Brieskorn modules} $\Htc_f^i$ by
$$\Htc_f^i:=\begin{cases}\Hc_f^i&\h{if}\,\,\,i\ne 1,\\
\Hc_f^{1}/\C\{t\}\omega_0&\h{if}\,\,\,i = 1.\end{cases}
\leqno(1.4.4)$$

Denote by $\Hc^i_{f,{\rm tors}}$ the $\dti$-torsion part of
$\Hc_f^i$, and set $\Hc^i_{f,{\rm free}}:=\Hc_f^i/\Hc^i_{f,{\rm tors}}$
(similarly for $\Htc^i_{f,{\rm tors}}$, $\Htc^i_{f,{\rm free}}$).
Let $R,K$ be as in (1.2.2). By [BaSa], Th.~1 we have
$$\aligned &\Hc^i_{f,{\rm tors}}=\Ker(\Hc_f^i\to\G_f^i),\\
&\dt^{-p}\,\Htc^i_{f,{\rm tors}}=0\q\h{if $\,p\gg 0\,$
(locally on $\,X_0$),}\\
&\Htc_{f,{\rm free},x}^i\cong\mopl^{\mu_{i,x}}R\q\h{with}\,\,\,
\mu_{i,x}:=\dim_K\GG_{f,x}^i.\endaligned
\leqno(1.4.5)$$

\ms\nin
{\bf 1.5.~Proof of Proposition~2.}
We may assume $c=0$ and moreover
$$\M=\OX,$$
since the assertion is analytic and local on $X$.
By (1.2.3) and (1.3.5) it is enough to show a canonical
isomorphism of sheaves of $\EE_{S,0}$-modules
$$\Hc^i\bl(\OX((\dti))^f\sotim_{\OX}\Omega_X^{\ssb}\br)\big|_{X_0}
=\Kh\sotim_K\GG_f^i,
\leqno(1.5.1)$$
in a compatible way with the canonical morphisms from $\G_f^i$.
\sk
With the notation of (1.3--4), consider the {\it algebraic}
microlocalization $\K_f^{\ssb}[\dti]$, and define
increasing filtrations $F_{\ssb}$ and $G_{\ssb}$ by
$$\aligned F_k\,\K_f^p[\dti]&:=\bl(\mopl_{i\les k+p}\,
\Omega_X^p\,\dt^i\br)\big|_{X_0},\\
G_k\,\K_f^p[\dti]&:=\A_f^p\,\dt^k\oplus\bl(\mopl_{i<k}\,
\Omega_X^p\,\dt^i\br)\big|_{X_0}.\endaligned
\leqno(1.5.2)$$
Note that $\A_f^{\ssb}$ is a subcomplex of $G_0\,\K_f^{\ssb}[\dti]$.
The {\it formal} microlocalization of $\K_f^{\ssb}$ can be defined by
$F$ and also by $G$, i.e.
$$\bl(\OX((\dti))^f\sotim_{\OX}\Omega_X^{\ssb}\br)\big|_{X_0}=
\rlap{\raise-8pt\h{$\,\,\scriptstyle p$}}\ilim
\bl(\rlap{\raise-8pt\h{$\,\,\,\scriptstyle q$}}\plim\,
(F_p/F_q)\,\K_f^{\ssb}[\dti]\br),
\leqno(1.5.3)$$
where $F$ can be replaced with $G$ since these are cofinal with each
other.
\sk
We have a canonical morphism between the short exact sequences
$$\begin{CD}0@>>>\A_f^{\ssb}@>>>\K_f^{\ssb}@>>>
\K_f^{\ssb}/\A_f^{\ssb}@>>>0\\
@. @VV{\al}V @VV{\beta}V @|\\
0@>>>G_0\,\K_f^{\ssb}[\dti]@>>>\K_f^{\ssb}[\dti]@>>>
\K_f^{\ssb}/\A_f^{\ssb}@>>>0\end{CD}$$
Using the morphism between the associated long exact sequences, we get
the canonical isomorphisms for any $i$
$$\Hc^i\bl(G_0\,\K_f^{\ssb}[\dti]\br)=\Htc_f^i,
\leqno(1.5.4)$$
where the $\Htc_f^i$ are as in (1.4.4).
Indeed, consider the long exact sequence associated to the first
short exact sequence, and then dividing $\Hc_f^1=\Hc^1\A_f^{\ssb}$ and
$\G_f^1=\Hc^1\K_f^{\ssb}$ by $\C[t]\omega_0$,
we still get a long exact sequence.
This is isomorphic to the long exact sequence associated to the
second short exact sequence by the five lemma.
Indeed, $\C[t]\omega_0\subset\G_f^1$ is the subsheaf annihilated by
the localization by $\dt$ by (1.3.5) and the action of $\dt$ on
$(\C\{t\}/\C[t])\omega_0$ is bijective so that the morphism induced
by $\Hc^1\beta$ is bijective.
\sk
We now calculate the cohomology sheaves of
$$\rlap{\raise-8pt\h{$\,\,\scriptstyle p$}}\ilim
\Bigl(\rlap{\raise-8pt\h{$\,\,\,\scriptstyle q$}}\plim\,
(G_p/G_q)\,\K_f^{\ssb}[\dti]\Bigr),$$
by using the Mittag-Leffler condition.
We first calculate the projective limit for $q$ with $p$ fixed.
Here we may assume $p=0$ by using the action $\dt$.
For $k'=k+i>k>0$, we have the commutative diagrams:
$$\begin{CD}0@>>>G_0\,\K_f^{\ssb}[\dti]@>{\dt^{-k'}}>>
G_0\,\K_f^{\ssb}[\dti]@>>>(G_0/G_{-k'})\,\K_f^{\ssb}[\dti]@>>>0\\
@. @VV{\dt^{-i}}V @| @VVV\\
0@>>>G_0\,\K_f^{\ssb}[\dti]@>{\dt^{-k}}>>G_0\,\K_f^{\ssb}[\dti]@>>>
(G_0/G_{-k})\,\K_f^{\ssb}[\dti]@>>>0\end{CD}$$
Combining this with (1.5.4) we get short exact sequences of
projective systems for $j\in\Z$
$$0\to\bl\{\Htc^j/\dt^{-k}\Htc^j\br\}_k\to
\bl\{\Hc^j\bl((G_0/G_{-k})\K_f^{\ssb}[\dti]\br)\br\}_k\to
\bl\{\Ker\bl(\dt^{-k}\big|\Htc^{j+1}\br)\br\}_k\to 0,$$
where the transition morphisms of the first and last projective
systems are induced by the identity on $\Htc^j$ and $\dt^{-(k'-k)}$
on $\Htc^{j+1}$ respectively.
By the middle property of (1.4.5), we see that the projective limit
of the last projective system vanishes. So we get isomorphisms for
$j\in\Z$
$$\rlap{\raise-9pt\h{$\,\,\,\scriptstyle k$}}\plim\,
\Hc^j\bl((G_0/G_{-k})\,\K_f^{\ssb}[\dti]\br)=\Hhc_f^j\q\h{with}\q
\Hhc_f^j:=\rlap{\raise-9pt\h{$\,\,\,\scriptstyle k$}}\plim\,
\Htc_f^j/\dt^{-k}\Htc_f^j.
\leqno(1.5.5)$$
Moreover, the Mittag-Leffler condition is satisfied for
$\bl\{\Hc^j\bl((G_0/G_{-k})\K_f^{\ssb}[\dti]\br)\br\}_k$
by using again the middle property of (1.4.5) together with the snake
lemma applied to the transition morphisms between the above short
exact sequence of projective systems.
So the projective limit commutes with the cohomology.
\sk
We have to calculate the inductive limit of the projective limits.
Here the inductive limit always commutes with the cohomology sheaf
functor.
The morphisms of the inductive system of the cohomology sheaves are
induced by the inclusions
$$G_p\,\K_f^{\ssb}[\dti]\into G_{p+i}\,\K_f^{\ssb}[\dti],$$
which are identified with the morphisms
$$\dt^{-i}:G_0\,\K_f^{\ssb}[\dti]\into G_0\,\K_f^{\ssb}[\dti],$$
under the identification
$$\dt^p:G_0\,\K_f^{\ssb}[\dti]\simto G_p\,\K_f^{\ssb}[\dti].$$
So (1.5.1) follows from (1.5.5) since the inductive limit of the
inductive system $\bl\{\Hhc_f^j\br\}_k$ with transition morphisms
defined by $\dt^{-(k'-k)}$ is the localization of $\Hhc_f^j$ by $\dti$.
This completes the proof of Proposition~2.

\ms\nin
{\bf 1.6.~Proof of Theorem~2.}
The assertion follows from Proposition~2 using (1.2.3) and (1.3.2).

\ms\nin
{\bf Remarks~1.7}~(i) The formal Brieskorn modules and Gauss-Manin
systems defined by the $\dti$-adic completion were used in
[Sai2] for hypersurface isolated singularities with non-degenerate
Newton polygons where the Mittag-Leffler condition was also used.
(See also Remark~(3.5) below.)
\ms
(ii) In case $f^{-1}(0)$ is a divisor with normal crossings, the
Brieskorn modules are $\dti$-torsion-free, and the filtration $G$ gives
the filtration $V$ of Kashiwara and Malgrange (indexed by $\Z$) by
using an argument similar to [St], (1.13), see also [Sai1].
It is possible to prove Proposition~2 in the normal crossing case
by using this together with (2.3) below.

\ms\nin
{\bf 1.8.~Commutativity of inductive limit and global cohomology.}
It is well known that the global cohomology functor does not necessarily
commute with inductive limit nor with $((u))$ in the noncompact case.
For instance, Cartan's Theorem B does not hold for quasi-coherent
sheaves, see e.g. [Sai3], 2.3.8.
So we give here a sufficient condition for their commutativity.
\sk
Let $\F_{\la}^{\ssb}$ be an inductive system of complexes of sheaves
on a topological space $X$ indexed by a directed set $\Lambda$.
(In this paper, inductive systems are always indexed by directed
sets.)
We have the canonical flasque resolutions
$$\F_{\la}^{\ssb}\to\I^{\ssb}\F_{\la}^{\ssb},\q
\F^{\ssb}\to\I^{\ssb}\F^{\ssb}\q\h{with}\,\,\,
\F^{\ssb}:=\rlap{\raise-9pt\h{$\,\,\scriptstyle\la$}}\ilim
\,\F_{\la}^{\ssb},$$
where $\I^{\ssb}$ denotes the canonical flasque resolution of Godement
using the discontinuous sections.
By the property of inductive limit, there is a canonical morphism
$$\rlap{\raise-9pt\h{$\,\,\scriptstyle\la$}}\ilim\,
\R\Gamma(X,\F_{\la}^{\ssb})\,(:=
\rlap{\raise-9pt\h{$\,\,\scriptstyle\la$}}\ilim\,
\Gamma(X,\I^{\ssb}\F_{\la}^{\ssb}))\to\R\Gamma(X,\F^{\ssb})\,
(:=\Gamma(X,\I^{\ssb}\F^{\ssb})).
\leqno(1.8.1)$$

Assume $X$ has a compactification $Y$ with a stratification
compatible with $X$ and satisfying the local triviality condition
using the local cone structure along each stratum.
More precisely, there is an increasing sequence of closed subspaces
$Y_k\,(k\ges -1)$ with $Y_{-1}=\emptyset$, $Y_d=Y\,\,(d\gg 0)$, and
for any $x\in Y_d\setminus Y_{d-1}$ with $d\ges 0$, there is an open
neighborhood $U_x$ of $x$ in $Y$ together with a compact topological
space $L_x$ having an increasing sequence of closed subspaces
$(L_x)_k\,(k\ges -1)$ with $(L_x)_{-1}=\emptyset$ and such that
there is a homeomorphism $U_x\cong\R^d\times C(L_x)$ inducing
$Y_k\cap U_x\cong\R^d\times C\bl((L_x)_{k-d-1}\br)$ for any
$k\ges d$, see e.g.\ [GoMa].
\sk
Set $Y_k^o:=Y_k\setminus Y_{k-1}$.
Assume furthermore the following condition is satisfied:
$$\h{$\Hc^i\F_{\la}^{\ssb}|_{Y_k^o}$ are locally constant for any
$\la,k,i$.}
\leqno(C)$$
We can then show that (1.8.1) is an quasi-isomorphism.
\sk
Let $j:X\into Y$ denote the inclusion.
We first prove the canonical quasi-isomorphism
$$\rlap{\raise-9pt\h{$\,\,\scriptstyle\la$}}\ilim\,
\R j_*\F_{\la}^{\ssb}\,(:=
\rlap{\raise-9pt\h{$\,\,\scriptstyle\la$}}\ilim\,
j_*\I^{\ssb}\F_{\la}^{\ssb})\simto
\R j_*\F^{\ssb}\,(:=j_*\I^{\ssb}\F^{\ssb}).
\leqno(1.8.2)$$
Since the inductive limit commutes with $\Hc^i$, we see that
condition~(C) holds also for $\F^{\ssb}$.
Using the local cone structure along each stratum, condition~(C)
then holds for $\R j_*\F_{\la}^{\ssb}$ and $\R j_*\F^{\ssb}$.
Consider the open inclusions
$$j_k:X\cup(Y\setminus Y_k)\into X\cup(Y\setminus Y_{k-1}).
\leqno(1.8.3)$$
Let $j'_k$ be the composition of $j_i$ for $i\ges k$.
Condition~$(C)$ for $\R(j'_{k+1})_*\F_{\la}^{\ssb}$ implies the
canonical isomorphisms for $x\in Y_k^0$
$$H^i(U_x\setminus Y_k,\R(j'_{k+1})_*\F_{\la}^{\ssb}|_{U_x\setminus
Y_k})\simto H^i(L_x,\R(j'_{k+1})_*\F_{\la}^{\ssb}|_{L_x}),$$
where $L_x$ is identified with $\{0\}\times\bl(\bl\{\frac{1}{2}
\br\}\times L_x\br)\subset\R^d\times C(L_x)\cong U_x$.
We have the same for $\R(j'_{k+1})_*\F^{\ssb}$.
So we can prove (1.8.2) at each stratum $Y_k^o$ by decreasing
induction on $k$ using the following well-known property:
$$\h{Inductive limit commutes with the global cohomology on compact
spaces.}
\leqno(1.8.4)$$

By (1.8.2) and (1.8.4), we get thus the canonical isomorphism
$$\rlap{\raise-9pt\h{$\,\,\scriptstyle\la$}}\ilim\,
\R\Gamma(X,\F_{\la}^{\ssb})\simto\R\Gamma(X,\F^{\ssb}).
\leqno(1.8.5)$$

\bs\bs
\centerline{\bf 2.~Formal microlocalization}
\bs\nin
In this section we recall some basics of formal microlocalization,
and prove Proposition~4 by reducing to the normal crossing case.

\ms\nin
{\bf 2.1.~Differential complexes and formal microlocalization.}
In this subsection $X$ is either a smooth complex algebraic variety or
a complex manifold.
By [Sai4], Prop.~1.8, there is an equivalence of categories
$$\aligned &\DR_X^{-1}:D^b(\OX,\Diff)\simto D^b(\D_X),\\
\h{with}\q&\DR_X^{-1}(\Lc):=\Lc\sotim_{\OX}\D_X\,\,\,
\h{for $\OX$-modules $\Lc$.}\endaligned$$
where $D^b(\D_X)$ is the bounded derived category of right
$\D_X$-modules, and $D^b(\OX,\Diff)$ is the derived category of
bounded complexes whose components are $\OX$-modules and whose
differentials are differential morphisms in the sense of loc.~cit.
More precisely, the differential morphisms are defined so that we have
for $\OX$-modules $\Lc,\Lc'$
$$\Hom_\Diff(\Lc,\Lc')=\Hom_{\D_X}(\DR_X^{-1}(\Lc),\DR_X^{-1}(\Lc')),$$
and the derived category is defined by inverting
$\D$-quasi-isomorphisms (which are, by definition, morphisms whose
mapping cones are $\D$-acyclic, where the last condition is defined by
the acyclicity of the image by the functor $\DR_X^{-1}$, see loc.~cit.)
Note that for $\LL\in D^b(\OX,\Diff)$ and $\M^{\ssb}\in D^b(\D_X)$,
there are canonical $\D$-quasi-isomorphisms and quasi-isomorphisms
$$\DR_X\bl(\DR_X^{-1}(\LL)\br)\simto\LL,\q
\DR_X^{-1}\bl(\DR_X(\M^{\ssb})\br)\simto\M^{\ssb},
\leqno(2.1.1)$$
where $\DR_X(\M)$ is the (shifted) de Rham complex for right
$\D_X$-modules $\M$.
\sk
Let $D^b_{rh}(\D_X)$ denote the full subcategory of $D^b(\D_X)$
consisting of $\M^{\ssb}$ with $\Hc^j\M^{\ssb}$ regular holonomic over
$\D_X$.
Let $D^b_{rh}(\OX,\Diff)$ be the full subcategory of
$D^b(\OX,\Diff)$ consisting of $\LL$ with
$\DR_X^{-1}(\LL)\in D^b_{rh}(\D_X)$.
We have an equivalence of categories
$$\DR_X^{-1}:D^b_{rh}(\OX,\Diff)\simto D^b_{rh}(\D_X).
\leqno(2.1.2)$$
By (2.1.1), any $\LL\in D^b_{rh}(\OX,\Diff)$ is represented by
$\DR_X\bl(\DR_X^{-1}(\LL)\br)$, and moreover, $\DR_X(\M)$ for
$\M\in M_{rh}(\D_X)$ are $\D$-quasi-isomorphic to subcomplexes
$F_p\DR_X(\M)$ (whose components are coherent over $\OX$) for
$p\gg 0$ locally on $X$ where $F$ is induced by a good filtration
on $\M$ locally defined on $X$.
By the functor $\DR^{-1}$, the direct image of $\D$-modules
is compatible with that of differential complexes, and the latter is
defined by the sheaf-theoretic direct image, see loc.~cit., Prop.~3.3.
\sk
Let $f$ be an algebraic or holomorphic function on $X$. Set $S:=\C$.
We have the graph embedding
$$i_f:X\into\XX:=X\times S.$$
We define the {\it formal microlocalization} $\LL((u))^f$ of
$\LL\in D^b_{rh}(\OX,\Diff)$ as follows.
\sk
Consider first $\DR_{\XX}^{-1}\LL$ where $\LL$ is identified with
the direct image by $i_f$.
It is a complex of right $\D_{\XX}$-module whose $j$th component is
$\Lc^j\sotim_{\OO_{\XX}}\D_{\XX}$.
Take the tensor product over $\D_{\XX}$ with
$$\D_{\XX\to S}:=\OO_{\XX}\sotim_{pr_2^{-1}\OS}pr_2^{-1}\D_S.$$
It has a structure of a complex of right $pr_2^{-1}\D_S$-modules.
We convert it to a complex of left $pr_2^{-1}\D_S$-modules by
using the involution $^*$ of $\D_S$ such that
$$(PQ)^*=Q^*P^*,\q g^*=g\,\,\,(g\in\OS),\q \dt^*=-\dt.
\leqno(2.1.3)$$
We thus get the first of the following three complexes:
$$\LL[\dt]^f,\q\LL[\dt,\dti]^f,\q\LL((\dti))^f.$$
The second complex is then obtained by taking the localization of
the first by $\dt$, and the third by the tensor product of the second
with $\OX((\dti))$ over $\OX[\dt,\dti]$ (which is well defined since
$\OX((\dti))$ is flat over $\OX[\dt,\dti]$).
Note that the above definition using the tensor product is reasonable
only in the case where the $\Hc^j\DR^{-1}(\LL)$ are coherent
(e.g. regular holonomic) over $\D_X$, see the remark after (2.1.2).
\sk
The upperscript $^f$ means that the differential and the action
of $t$ are twisted by using $f$ as in the introduction, and
$\LL((\dti))^f$ may be understood as an abbreviation of
$\LL((\dti))\delta(f-t)$ where $\delta(f-t)$ is the delta
function satisfying the relations
$$t\,\delta(f-t)=f\,\delta(f-t),\q\xi\,\delta(f-t)=-(\xi f)\,
\delta(f-t)\q\h{for}\,\,\xi\in\Theta_X.$$

In case $\LL=\DR_X(\M)$ for a locally finite free $\OX$-module $\M$
with an integrable connection, we have a canonical isomorphism
$$\LL((\dti))^f=\DR_X\bl(\M((\dti))^f\br).
\leqno(2.1.4)$$
Moreover, the construction of $\LL((\dti))^f$ is compatible
with the direct image by a proper morphism $\pi:X'\to X$, i.e.
we have for a differential complex $\LL$ on $X'$
$$\R\pi_*\bl(\LL((\dti))^{\pi^*f}\br)=
(\R\pi_*\LL)((\dti))^f,
\leqno(2.1.5)$$
where $\R\pi_*$ is defined by the sheaf-theoretic direct image.
Here we may assume that each component of the complexes is
$\pi_*$-acyclic by taking the canonical flasque resolution of
Godement.
It is then enough to show the compatibility with the differential
of the complexes where the formal microlocalization can be replaced
with the algebraic one by the above argument.
So the assertion is shown by using the direct image of induced
$\D$-modules as in the proof of [Sai4], Prop.~3.3.
\sk
We will also consider the case where the divisor at infinity $D$ is
given in $\X$ (especially in the analytic case).
In this case, we similarly define $D^b_{rh}(\OXb(*D),\Diff)$ by
assuming that the $\Lb{}^j$ are $\OXb(*D)$-modules and the
$\Hc^j\DR_X^{-1}(\Lb{}^{\ssb})$ are regular holonomic over $\D_{\X}$,
and then define the formal microlocalization $\Lb{}^{\ssb}((\dti))^f$
by using the tensor product with $\OXb(*D)((\dti))^f$ over $\OXb(*D)$
instead of the one with $\OXb((\dti))^f$ over $\OXb$.
(For the difference between these, see Example~(2.2) below.)
We have also the compatibility with the direct image as above by
setting $D':=\bar{\pi}{}^{-1}(D)$ for a proper morphism
$\bar{\pi}:\X'\to\X$.
These will be used at the end of (2.3) below.

\ms\nin
{\bf Example~2.2.} Set $X:=\C$, $E:=\{0\}$, and also $\X:=\C$,
$D:=\{0\}$. Let $x$ be the coordinate of $\C$. Set
$$\LL:=C(x\dd_x:\OX\to\OX),\q
\Lb{}^{\ssb}:=C(x\dd_x:\OXb(*D)\to\OXb(*D)).$$
We have
$$\D_X/\D_X\,\dd_xx=\OX[x^{-1}]=\OX(*E),$$
and there are canonical quasi-isomorphisms
$$\DR_X^{-1}(\LL)\simto\OX(*E),\q
\DR_{\X}^{-1}(\Lb{}^{\ssb})\simto\OXb(*D),$$
where the last one follows from the first by using the localization.
(The difference between $x\dd_x$ and $\dd_xx$ comes from the involution
$^*$ of $\D_X$ in (2.1.3) for the transformation between left and right
$\D$-modules.)
\sk
However, setting $f=x$, we can easily show by a direct calculation
$$\LL((\dti))^f\simto\C((\dti))_0,\q\Lb{}^{\ssb}((\dti))^f=0,$$
where $\C((\dti))_0$ is the sheaf supported at the origin with stalk
$\C((\dti))$.

\ms\nin
{\bf 2.3.~Reduction of Proposition~4 to the normal crossing case.}
Let $\pib:\X'\to\X$ be a proper birational morphism of smooth
complex algebraic varieties such that the union of the pullback $D'$
of $D$ and a finite number of fibers $\X'_c\,\,(c\in\Sigma)$ is a
divisor with simple normal crossings (where $\Sigma\subset\PP^1$ is a
finite set containing $\infty$), and moreover the restriction of
$(\X',D')$ over $S':=S\setminus\Sigma$ is a divisor with simple
normal crossings over $S'$.
(The latter condition means that the restriction over $S'$ of any
intersection of irreducible components of $D'$ is smooth over $S'$.)
Set $X':=\pib^{-1}(X)$, $\pi:=\pib|_{X'}:X'\to X$, and
$\M':=\pi^*\M$.
There is a canonical morphism of differential complexes
$$\pi^{\#}:\DR_X(\M)\to\R\pi_*\DR_{X'}(\M').
\leqno(2.3.1)$$
Applying this to the dual $\M^*:=\Hom_{\OX}(\M,\OX)$, and using
the duality for the direct images of differential complexes
(see [Sai4], Th.~3.11), we get 
$$\pi_{\#}:\R\pi_*\DR_{X'}(\M')\to\DR_X(\M).
\leqno(2.3.2)$$
Since the composition $\pi_{\#}\ssc\pi^{\#}$ is the identity on
$\DR_X(\M)$, we get an isomorphism
$$\R\pi_*\DR_{X'}(\M')\cong\DR_X(\M)\oplus{\rm Cone}\,\pi^{\#}\q
\h{in}\,\,D^b(\OX,\Diff)^f,
\leqno(2.3.3)$$
together with
$${\rm Cone}\,\pi^{\#}\cong{\rm Cone}\,\pi_{\#}[-1].$$
Indeed, these follow from the octahedral axiom of the triangulated
category (applied to the composition of $\pi^{\#}$ and $\pi_{\#}$).
\sk
Let $j':X'\into\X'$ denote the inclusion. Set $\MM{}':=j_*\M'$.
By (2.3.3) and GAGA, we have
$$\DR_{\X^{\an}}(\MM{}^{\an})\,\,\,\h{is a direct factor of}
\,\,\,\R\pib_*\DR_{\X'{}^{\an}}(\MM{}'{}^{\an}).
\leqno(2.3.4)$$
Set $f':=\pi^*f$. By (2.1.5) together with the remark after it,
we then get
$$\DR_{\X^{\an}}(\MM{}^{\an})((\dti))^f\,\,\,
\h{is a direct factor of}\,\,\,\R\pib_*\bl(\DR_{\X'{}^{\an}}
(\MM{}'{}^{\an})((\dti))^{f'}\br).
\leqno(2.3.5)$$
Here we use the remark at the end of (2.1).
Thus Proposition~4 is reduced to the case where the union of
$D:=\X\setminus X$ and the fibers $\X_c\,\,(c\in\Sigma)$ is a divisor
with normal crossings.

\ms\nin
{\bf 2.4.~Proof of Proposition~4 in the normal crossing case.}
Since the assertion is local on $\X$, we may restrict to an open
neighborhood in $\X$ of a point of $\X_c\setminus X_c$,
where we may assume $c=0$ or $\infty$ (by replacing $f$ with $f-c$
in case $c\ne\infty$).
By (2.3) we may assume that $X$, $\X$ are complex manifolds of
dimension $n$ such that the union of $D:=\X\setminus X$ and
$f^{-1}(0)$ is a divisor with normal crossings on $\X$.
We take a neighborhood of $0\in\X_c$ in $\X$ of the form
$$U_{\ep}=\De_{\ep}^n\subset\X,$$
where $\De_{\ep}$ is an open disk of radius $\ep$.
We may assume that there are subsets $J,J'\subset\{1,\dots,n\}$
such that
$$U'_{\ep}:=U_{\ep}\cap X=\bl\{\mprod_{i\in J'}\,x_i\ne 0\br\},\q
f=\mprod_{i\in J}\,\1x_i^{e_i},$$
replacing $U_{\ep}$ and $\ep$ if necessary, where $x_1,\dots,x_n$ are
the coordinates of the polydisk $\De_{\ep}^n$. Note that
$e_i\les-1$ and $J'\supset J$ if $c=\infty$, and $e_i\ges 1$ if $c=0$.
In the inductive argument below, we have to treat also the
case $f=0$ where we set $J=\emptyset$.
We assume $J\ne\emptyset$ (i.e. $f\ne 0$) for the moment until the
proof for the case $J=\emptyset$ will be explained.
\sk
We first treat the case $c=0$.
Since the assertion is local on $\X$, we may replace $\X$, $X$ with
$U_{\ep}$, $U'_{\ep}$ respectively, and we will denote by $\MM$, $\M$
their restrictions to $U_{\ep}$, $U'_{\ep}$ in this subsection.
The assertion is reduced to the case ${\rm rank}\,\M=1$ by using a
filtration on $\M$ since the local monodromies are commutative in
the normal crossing case. Then there is a generator $m$ of $\MM$
annihilated by the differential operators
$$\xi_i:=\begin{cases}\dd_{x_i}x_i-\al_i&\h{if}\,\,\,i\in J',\\
\dd_{x_i}&\h{if}\,\,\,i\notin J',\end{cases}$$
where $\al_i\in\C$, and we have $\MM=\D_{U_{\ep}}/\J$ with $\J$
generated by the above differential operators.
\sk
Replacing $\xi_i=\dd_{x_i}$ with $\xi'_i:=\dd_{x_i}x_i$ for
$i\in J\setminus J'$, we get a left ideal $\J'$ of $\D_{U_{\ep}}$
generated by the $\xi'_i$, where $\xi'_i:=\xi_i$ if
$i\notin J\setminus J'$. Then there is an injection
$$\MM\into\MM{}':=\D_{U_{\ep}}/\J',$$
defined by the multiplication by $\prod_{i\in J\setminus J'}x_i$.
Moreover, there is a finite increasing filtration $W$ on $\MM{}'$
such that
$$\Gr^W_k\MM{}'=\mopl_{I\subset J\setminus J',\,|I|=k}\,\MM_I,$$
where the $\MM_I$ are $\D_{U_{\ep}}$-modules of the same type as $\MM$
(up to the direct image as $\D$-modules by closed immersions), and
are supported on
$$U_{\ep,I}:=\mcap_{i\in I}\,\{x_i=0\}\subset U_{\ep}.$$
If $|I|=k>0$, then $U_{\ep,I}$ is contained in $f^{-1}(0)$.
The assertion in this case is reduced to the case $f=0$ (i.e.
$J=\emptyset$) by using differential complexes as in (2.1), and can
be proved by using the inductive argument in this subsection.
We may thus replace $\MM$ with $\MM{}'$.
\sk
Since the $\xi'_i$ give a free resolution of $\MM{}'$,
we get an isomorphism
$$\DR_{U_{\ep}}(\MM{}')=\Kb(\OO_{U_{\ep}};\xi_i^{\prime *})[n],$$
where the right-hand side is the (shifted) Koszul complex associated
with the mutually commuting operators $\xi_i^{\prime *}$ with $^*$
the involution of $\D_{U_{\ep}}$ as in (2.1.3). (Indeed, the
right-hand side is actually $\Kb(\Omega_{U_{\ep}}^n;\xi'_i)[n]$.)
Here we may replace $\xi_i^{\prime *}$ with $-\xi_i^{\prime *}$.
\sk
Let $g:=\prod_{i\in J'}x_i$ so that
$U_{\ep}\setminus U'_{\ep}=g^{-1}(0)$. Set
$$A_{\ep}:=\Gamma(U_{\ep},\OO_{U_{\ep}})
\bl[\h{$\frac{1}{g}$}\br],\q
A'_{\ep}:=\Gamma(U'_{\ep},\OO_{U'_{\ep}})\q
A''_{\ep}=A'_{\ep}/A_{\ep}.$$
By the above calculation of the de Rham complex $\DR_{U_{\ep}}(\MM{}')$
using the Koszul complex of the operators $-\xi_i^{\prime *}$,
its twisted de Rham complex is given by the Koszul complex of the
following operators:
$$\xit_i:=\begin{cases}
x_i\dd_{x_i}+\al_i+e_if\dt&\h{if}\,\,\,i\in J\cap J',\\
x_i\dd_{x_i}+e_if\dt&\h{if}\,\,\,i\in J\setminus J',\\
x_i\dd_{x_i}+\al_i&\h{if}\,\,\,i\in J'\setminus J,\\
\dd_{x_i}&\h{if}\,\,\,i\notin J\cup J'.\end{cases}$$
More precisely, let
$$\Kb\bl(A_{\ep}((\dti)),\xit_i\br),\q
\Kb\bl(A'_{\ep}((\dti)),\xit_i\br),\q
\Kb\bl(A''_{\ep}((\dti)),\xit_i\br)$$
respectively denote the Koszul complex associated with the mutually
commuting differential operators $\xit_i$ on
$A_{\ep}((\dti))$, $A'_{\ep}((\dti))$ and $A''_{\ep}((\dti))$.
For the proof of Proposition~4, we have to show the canonical
quasi-isomorphism
$$\MM((\dti))^f\sotim_{\OO}\Omega_{U_{\ep}}^{\ssb}\to
\R j_*\bl(\M((\dti))^f\sotim_{\OO}\Omega_{U'_{\ep}}^{\ssb}\br).$$
We first show that the stalk at the origin of this morphism is
identified with the inductive limit for $\ep\to 0$ of the canonical
morphism
$$\Kb\bl(A_{\ep}((\dti)),\xit_i\br)\to
\Kb\bl(A'_{\ep}((\dti)),\xit_i\br).
\leqno(2.4.1)$$
Indeed, the assertion for $\Kb\bl(A_{\ep}((\dti)),\xit_i\br)$ follows
from the {\it definition} of the projective limit of sheaves.
As for $\Kb\bl(A'_{\ep}((\dti)),\xit_i\br)$, set
$$\Kch_f^{\1\prime\ssb}:=\M'((\dti))^f\sotim_{\OO}
\Omega^{\ssb}_{U'_{\ep}}\q\h{with}\q\M':=\MM{}'|_{U'_{\ep}},$$
where the differential is defined as in the introduction.
We have isomorphisms of complexes
$$\Kb\bl(A'_{\ep}((\dti)),\xit_i\br)=
\Gamma\bl(U'_{\ep},\Kch_f^{\1\prime\ssb}\br)=
\rlap{\raise-8pt\h{$\,\,\scriptstyle p$}}\ilim\,
\Gamma\bl(U'_{\ep},G_p\,\Kch_f^{\1\prime\ssb}\br),
\leqno(2.4.2)$$
where the filtration $G$ on $\Kch_f^{\1\prime\ssb}$ is defined in the
same way as in (1.5). We have moreover a canonical quasi-isomorphism
$$\Gamma\bl(U'_{\ep},G_p\,\Kch_f^{\1\prime\ssb}\br)\simto
\R\Gamma\bl(U'_{\ep},G_p\,\Kch_f^{\1\prime\ssb}\br),
\leqno(2.4.3)$$
since the infinite direct product of sheaves $\prod$ commutes with
cohomology, and $U'_{\ep}$ is Stein.
We can show that their restrictions to the strata of the natural
stratification of a divisor with normal crossings are locally constant
by reducing to the case of Remark~(1.7)(ii) using the filtration $W$
on $\MM{}'$ defined above. So we can apply (1.8) to get the
commutativity of the inductive limit and $\R\Gamma(U'_{\ep},*)$.
Combining this with (2.4.2--3), we get the quasi-isomorphism
$$\Kb\bl(A'_{\ep}((\dti)),\xit_i\br)\simto
\R\Gamma\bl(U'_{\ep},\Kch_f^{\1\prime\ssb}\br).$$
The proof of Proposition~4 in the normal crossing case is thus reduced
to the assertion that (2.4.1) is a quasi-isomorphism.
Note that this is equivalent to
$$\Kb\bl(A''_{\ep}((\dti)),\xit_i\br)=0.
\leqno(2.4.4)$$

Here we may assume
$$J\setminus J'\ne\emptyset,\,\,\,\h{i.e.}\,\,\,J\not\subset J'.
\leqno(2.4.5)$$
Indeed, if $J\subset J'$, then the target of (2.4.1) vanishes, and
we get the vanishing of the source by using the $\dti$-adic
filtration on $A_{\ep}((\dti))$ defined by $A_{\ep}[[\dti]]\dt^{-k}$
for $k\in\Z$. Here the graded piece of $\xit_i$ for $i\in J$ is
given by $e_if\dt$ and $f$ is invertible in $A_{\ep}$ by the condition
$J\subset J'$. (A similar argument applies to the case $c=\infty$.)
\sk
For an element of $A_{\ep}((\dti))$ or $A'_{\ep}((\dti))$, we have
the Laurent expansion
$$\msum_{\nu\in\Z^n}\,a_{\nu}x^{\nu}\q\h{with}\,\,\,a_{\nu}\in\C,$$
where the coefficients $a_{\nu}$ satisfy certain convergence
conditions as is well-known. Since
$$(x_i\dd_{x_i}-\nu_i)\,x^{\nu}=0,$$
we see that $\Kb\bl(A_{\ep}((\dti)),\xit_i\br)$,
$\Kb\bl(A'_{\ep}((\dti)),\xit_i\br)$ are acyclic in case
$\al_i\notin\Z$ for some $i\in J'\setminus J$ (using the $n$-ple
complex structure of the Koszul complex). We may thus assume
$$\al_i\in\Z\,\,\,\h{for any}\,\,\,i\in J'\setminus J.$$
By the above calculation, the assertion is further reduced to the
case $J'\setminus J=\emptyset$ by taking the kernel and cokernel of
$\xit_i$ for $i\in J'\setminus J$ inductively (and using the $n$-ple
complex structure of the Koszul complex). Here we also use the external
product $\OO_{X_1\mtim X_2}=\OO_{X_1}\boxtimes\OO_{X_2}$ together with
Theorem~2 in case $J'\cap J=\emptyset$. We may thus assume
$$J'\setminus J=\emptyset,\,\,\,\h{i.e.}\,\,\,J'\subset J.$$
By a similar argument we may assume moreover
$$[1,n]\setminus J=\emptyset,\,\,\,\h{i.e.}\,\,\,J=[1,n].$$
Note that $n=|J|\ges 2$, since $J\setminus J'\ne\emptyset$ and
$J'\ne\emptyset$.
\sk
Choose $i_0\in J\setminus J'$. For the calculation of the Koszul
complex, we may replace $\xit_i$ with $\xit'_i$ defined by
$$\xit'_i:=\begin{cases}\xit_i&\h{if}\,\,\,i=i_0,\\
\xit_i-(e_i/e_{i_0})\,\xit_{i_0}&\h{if}\,\,\,i\ne i_0,\end{cases}$$
(since this does not change the vector space spanned by the $\xit_i$.)
For $i\ne i_0$, we have
$$\xit'_i=x_i\dd_{x_i}-(e_i/e_{i_0})\,x_{i_0}\dd_{x_{i_0}}-\al'_i
\q\h{with}\,\,\,\al'_i\in\C,$$
and hence
$$\xit'_i\,x^{\nu}=\bl(\nu_i-(e_i/e_{i_0})\,\nu_{i_0}-\al'_i\br)\,
x^{\nu}.
\leqno(2.4.6)$$

For $I\subset[1,n]\setminus\{i_0\}$, define
$$A_{\ep,I}:=\bl\{\msum_{\nu\in\Z^n}\,a_{\nu}x^{\nu}\in A_{\ep}
\,\,\big|\,\,a_{\nu}=0\,\,\,\h{if}\,\,\,\nu_i\ne(e_i/e_{i_0})\,\nu_{i_0}
+\al'_i\,\,\,\h{for some}\,\,\,i\in I\br\},$$
and similarly for $A'_{\ep,I}$ with $A_{\ep}$ replaced by
$A'_{\ep}$. Consider a canonical morphism
$$\Kb\bl(A_{\ep,I}((\dti)),\xit_i\br)\to
\Kb\bl(A'_{\ep,I}((\dti)),\xit_i\br).$$
By using (2.4.6) together with the $n$-ple complex structure of the
Koszul complex, the assertion (2.4.4) is then reduced by increasing
induction on $|I|$ to the assertion that the above canonical morphism
is a quasi-isomorphism for some $I\subset[1,n]\setminus\{i_0\}$.
But in case $I=[1,n]\setminus\{i_0\}$, we have the isomorphism
$$A_{\ep,I}\simto A'_{\ep,I}.$$
(Indeed, this easily follows from the condition: $i_0\in J\setminus J'$.)
So the assertion is proved for $c=0$, except for the case $J=\emptyset$
(i.e. $f=0$) which is also needed in our inductive argument. However,
the proof in the case $J=\emptyset$ is similar, and is easier since
the replacement of $\xit_i$ by $\xit'_i$ is not needed in this case.
(The details are left to the reader.)
\sk
The assertion for $c=\infty$ is equivalent to the vanishing of the
restriction of the left-hand side of the isomorphism in Proposition~4
to $\X_{\infty}$. We have $J\subset J'$ in this case. So the assertion
is proved in the same way as in the explanation after (2.4.5) by using
the $\dti$-adic filtration. This finishes the proof of Proposition~4.

\bs\bs
\centerline{\bf 3. Proof of the main theorem.}

\bs\nin
In this section we prove Theorem~1 after showing Propositions~1 and 3.

\ms\nin
{\bf 3.1.~Proof of Proposition~1.}
Using the continuous morphism $\rho:\X^{\an}\to\X$,
we get the canonical morphisms
$$\DR_{\X}\bl(\MM((u))^f\br)\buildrel{\rho^{\#}}\over\to
\rho_*\DR_{\X^{\an}}\bl(\MM^{\an}((u))^f\br)\to\R
\rho_*\DR_{\X^{\an}}\bl(\MM^{\an}((u))^f\br).
\leqno(3.1.1)$$
The middle term is defined by using
$\Gamma\bl(U^{\an},\Omega_{U^{\an}}^p\sotim_{\OO_{U^{\an}}}
\MM^{\an}((u))^f|_{U^{\an}}\br)$ for sufficiently small affine
open subvarieties $U$ of $\X$.
Here the non-commutativity of inductive limit and cohomology
does not cause a problem for the construction of the morphism
$\rho^{\#}$.
\sk
Taking the global sections on $\X$, we get the canonical morphism
in Proposition~1.
We can then use the truncation $\sigma_{\ges p}$, and the assertion
is reduced to the case of a sheaf $\F((u))$ with
$$\F=\MM\sotim_{\OXb}\Omega_{\X}^p\q(p\in\Z).$$
We have the decomposition
$$\F((u))=\F[[u]]\oplus\F[u^{-1}]u^{-1},
\leqno(3.1.2)$$
and the global cohomology functors on $\X$ and $\X^{\an}$ commute with
the infinite direct product and also with the infinite direct sum.
So the assertion is reduced to the case of $\F$, and follows from
GAGA.
This finishes the proof of Proposition~1.

\ms\nin
{\bf 3.2.~Proof of Proposition~3.} We may assume $c=0$. Set
$$\K_f^{\ssb}:=\bl(\M[\dt]^f\sotim_{\OX}\Omega_X^{\ssb}\br)\big|_{X_0}.
\leqno(3.2.1)$$
This coincides with the notation in (1.3) if $\M=\OX$.
We have the filtration $V$ on $\K_f^{\ssb}$ induced by the
filtration $V$ of Kashiwara and Malgrange on $\M[\dt]^f$, and
the assertion is equivalent to the following condition:
$$\h{The induced filtration $V$ on
$\R\Gamma\bl(X_0,\K_f^{\ssb}\br)$ is strict.}
\leqno(3.2.2)$$
We have a spectral sequence of regular holonomic $\D_{S,0}$-modules
$$E_2^{p,q}=H^p\bl(X_0,\Hc^q\K_f^{\ssb}\br)
\Rightarrow H^{p+q}\bl(X_0,\K_f^{\ssb}\br).
\leqno(3.2.3)$$
This is a filtered spectral sequence since $\K_f^{\ssb}$ has
the filtration $V$ induced by the filtration $V$ of Kashiwara and
Malgrange on $\M[\dt]^f$.
More precisely, let $\E$ denote the category of $\C$-vector spaces
having an exhaustive increasing filtration indexed by $\Z$,
and $\A$ be the abelian category consisting of graded vector spaces
$\mopl_{i\in\Z}\,E_i$ endowed with morphisms $u_i:E_i\to E_{i+1}$.
Then $\E$ is identified with the full subcategory of $\A$ consisting
of subobjects such that the $u_i$ are all injective.
The spectral sequence (3.2.3) is defined in $\A$, since $V$ is
essentially indexed by $m^{-1}\Z\subset\Q$ for some nonzero integer
$m$.
\sk
By [De2], 1.3 and [Sai1], 1.3.6, the assertion of Proposition~3
is then reduced to
$$E_r^{p,q}\in\E\q\h{for any $p,q\in\Z$ and $r\ges 2$}.
\leqno(3.2.4)$$
We show this by increasing induction on $r\ges 2$.
Assume $r=2$.
The induced filtration on the cohomology sheaves $\Hc^q\K_f^{\ssb}$
is strictly induced, see (1.3.3).
It is also strictly induced on
$E_2^{p,q}=H^p\bl(X_0,\Hc^q\K_f^{\ssb}\br)$ by taking the
inductive limit for $j\to-\infty$ of the canonical splitting (1.1.7)
which also holds for constructible sheaves of $\D_{S,0}$-modules.
So the assertion is shown for $r=2$.
\sk
We can then proceed by increasing induction on $r\ges 2$ using
the property that the $V$-filtration is strictly compatible with
any morphism of regular holonomic $\D_{S,0}$-modules,
see e.g.\ [Sai1], 3.1.5.
This finishes the proof of Proposition~3.

\ms\nin
{\bf 3.3.~Proof of Theorem~1.}
By Propositions~1 and 4, we may replace $X$ and $\M$ with the
associated analytic objects, and we will denote them by $X$
and $\M$ without adding $^{\an}$.
In particular, we consider only analytic sheaves.
\sk
We have the filtration $V$ of Kashiwara [Kas1] and Malgrange [Mal]
on $\M[\dt]^f$ such that the roots of the minimal polynomial of the
action of $t\dt$ on $\Gr_V^i\M[\dt]^f$ are contained in $\La+i$
where $\La$ is as in (1.1.4). Set $n=\dim X$.
Then we have the canonical isomorphisms
$$\psi_fL[n-1]=\DR_X(\Gr_V^1\M[\dt]^f),\q
\varphi_fL[n-1]=\DR_X(\Gr_V^0\M[\dt]^f),
\leqno(3.3.1)$$
such that the monodromy $T$ on the left-hand side corresponds to
$\exp(-2\pi it\dt)$ on the right-hand side, see loc.~cit.
Here $\DR_X$ is shifted by $n$ to get perverse sheaves as usual.
By Proposition~2 and (1.2.3), the assertion is then reduced to the
canonical isomorphisms
$$\EE_{S,0}\sotim_{\D_{S,0}}H^k\bl(X_0,\K_f^{\ssb}\br)\simto
H^k\bl(X_0,\EE_{S,0}\sotim_{\D_{S,0}}\K_f^{\ssb}\br),
\leqno(3.3.2)$$
where $\K_f^{\ssb}:=\bl(\M[\dt]^f\sotim_{\OX}\Omega_X^{\ssb}\br)
|_{X_0}$ as in (3.2.1).
Here we may assume $c=0$ as in (3.1).
\sk
To show (3.3.2), consider the spectral sequence
$$\Eh_2^{p,q}=H^p\bl(X_0,\Hc^q(\EE_{S,0}\sotim_{\D_{S,0}}
\K_f^{\ssb})\br)\Rightarrow
H^{p+q}\bl(X_0,\EE_{S,0}\sotim_{\D_{S,0}}\K_f^{\ssb}\br).
\leqno(3.3.3)$$
We have a canonical morphism form the spectral sequence (3.2.3)
to (3.3.3).
So it is enough to show that it induces the isomorphisms
$$\EE_{S,0}\sotim_{\D_{S,0}}E_r^{p,q}\simto\Eh_r^{p,q}.
\leqno(3.3.4)$$
For $r=2$, the isomorphism follows from (1.2.3).
Then we can proceed by increasing induction on $r\ges 2$ since
$\EE_{S,0}$ is flat over $\D_{S,0}$.
So (3.3.4) and then (3.3.3) follow.
This finishes the proof of Theorem~1.

\ms\nin
{\bf 3.4.~Cech calculation.}
Let $\{U_i\}$ be an affine open covering of $X$.
Set $U_I:=\mcap_{i\in I}\,U_i$ for $I\ne\emptyset$.
Then the left-hand side of the formula in Theorem~1 is given by
the $k$-th cohomology of the single complex associated with a
double complex whose $(p,q)$-component is
$$\mopl_{|I|=p+1}\,\Gamma\bl(U_I,\M\sotim_{\OO_X}\Omega_X^q|_{U_I}\br)
((u))^f,
\leqno(3.4.1)$$
where the first and second differentials are respectively given by
the Cech differential and $\nabla-u^{-1}df\wedge$.

\ms\nin
{\bf 3.5.~Isolated singularity case.} If $f$ has only isolated
singular points, then the filtration $F$ is strict, i.e.\
we have the $E_1$-degeneration of the spectral sequence associated
with the filtration $F$ on the left-hand side of the formulas in
Theorems~1 and 2, where $F$ is defined by a natural generalization of
(1.5.2) to the local system case. In fact, their $E_0$-complexes are
essentially given by the Koszul complexes associated with $df\wedge$
which are acyclic except for the top degree.
This implies that we get the $\dti$-adic completion of the Gauss-Manin
system at each singular point of $f$, and the filtration $F$ on it
is given by applying $\dt^k$ ($k\in\Z$) to the completion of the
Brieskorn lattice up to a shift by $\dim X$.
(See also Remark~(1.7)(i) for the relation with [Sai2].)

\ms\nin
{\bf Remark~3.6.} Let $V^k_c$, $V^k$, $A_1(c,c)$ be as in the remark after
Theorem~1. By $[t,\dt^{-i}]=i\,\dt^{-i-1}$, we can modify $V^k$ with
$\Rh\,V^k$ and $A_1(c,c)$ unchanged, so that we have by
increasing induction on $i\ges 1$
$$(t-c)\,V^k_c\subset \dti\Rh\,V^k_c+\dt^{-i-1}\Rh\,V^k.$$
In fact, if $\{v_{c,j}\}$ is a basis of $V^k_c$, then we can replace
$v_{c,j}$ with
$$v_{c,j}+\msum_{j',c'\ne c}\,a_{c,j,c',j'}\,\dt^{-i}\,v_{c',j'},$$
by increasing induction on $i\ges 1$ where the $a_{c,j,c',j'}\in\C$
are appropriately chosen. 
Passing to the limit, this means that $A_i(c,c')=0$ for $c\ne c'$.
The conclusion of the remark after Theorem~1 then follows.

\bs\bs
\centerline{\bf 4.~Generalization to the regular holonomic case}
\bs\nin
In this section we show how to generalize the arguments in the
previous sections in order to prove Theorem~3.

\ms\nin
{\bf 4.1.~Gauss-Manin systems with coefficients.}
Let $X$ be a complex manifold, and $E$ be a divisor with simple
normal crossings on it.
Let $\M$ be a locally free $\OX$-module of finite type endowed
with an integrable regular singular connection $\nabla$ having simple
poles along $E$ such that the real part of any eigenvalue of the
residue of $\nabla$ along each irreducible component of $E$ is
contained in $[0,1)$.
Let $\M(*E)$ be the localization of $\M$ along $E$.
This is a regular holonomic left $\D_X$-module.
Let $L$ be the local system on $X\setminus E$ corresponding to
$\M|_{X\setminus E}$. We have the canonical isomorphisms in
$D^b_c(X,\C)$:
$$\R j'_*L\simto\M\sotim_{\OX}\Omega_X^{\ssb}(\log E)\simto
\M(*E)\sotim_{\OX}\Omega_X^{\ssb},$$
where $j':X\setminus E\into X$ is the inclusion, see [De1].
\sk
Let $f$ be a holomorphic function on $X$ with $X_0:=f^{-1}(0)\subset E$.
The direct image as a left $\D$-module of $\M(*E)$ by the graph
embedding of $f$ is identified with
$$\M(*E)[\dt]^f,$$
where $^f$ means that the action of $\D_X$ and $t$ are twisted as in
(1.3). Set
$$\K_{\M,f}^{\ssb}:=\M[\dt]^f\sotim_{\OX}\Omega_X^{\ssb}(\log E)|_{X_0}
\subset\M(*E)[\dt]^f\sotim_{\OX}\Omega_X^{\ssb}|_{X_0}.$$
Here the last complex is the relative de Rham complex associated with
$\M(*E)[\dt]^f$, and the last inclusion is a quasi-isomorphism.
The Gauss-Manin systems with coefficients in $\M$ are defined by
$$\G_{\M,f}^i:=\Hc^i\K_{\M,f}.$$
These are constructible sheaves of regular holonomic $\D_{S,0}$-modules
on $X_0$ as in (1.3), and
$$\bl(\G^{i+1}_{\M,f}\br)^{\al}=\Hc^i\psi_{f,\la}\R j'_*L\,\,\,\h{with}
\,\,\,\la:=\exp(-2\pi i\al),
\leqno(4.1.1)$$
where the left-hand side is defined by
$\Ker(t\dt-\al)^k\subset\G^{i+1}_{\M,f}$ for $k\gg 0$ as in (1.1).
Note that $\psi_{f,\la}\R j'_*L=\varphi_{f,\la}\R j'_*L$ in this case.
Moreover, the filtration $V$ on the $\G_{\M,f}^i$ is induced by the
filtration $V$ on the $\D$-module $\M(*E)[\dt]^f$. These can be shown by
using the Milnor fibration around each point of $X_0$, see e.g. [Sai3],
Prop.~3.4.8.
\sk
Define a subcomplex $\A_{\M,f}^{\ssb}\subset\K_{\M,f}^{\ssb}$ by
$$\A_{\M,f}^i:=\Ker(df\wedge:\M\sotim_{\OX}\Omega_X^i(\log E)\to
\M\sotim_{\OX}\Omega_X^{i+1}(\log E))|_{X_0}.$$
We show that $\Hc^i\bl(\A_{\M,f}^{\ssb}\br)_x$ for $x\in X_0$
is a finite free $\C\{t\}$-submodule of $\bl(\G_{\M,f}^i\br)_x$ which
generates it over $\D_{S,0}$ and is stable by the action of $t\dt$
so that the real part of any eigenvalue of the residue of $t\dt$ is
contained in the $[-1,0)$. Setting $\la:=\exp(-2\pi i\al)$, this
assertion is equivalent by (4.1.1) to the following:
$$\bl(\Hc^{i+1}\A_{\M,f}\br)^{\al}=\begin{cases}
\Hc^i\psi_{f,\la}\R j'_*L&\h{if}\,\,\,\al\ges01,\\
0&\h{if}\,\,\,\al<-1.\end{cases}
\leqno(4.1.2)$$
Using this, we can define the filtration $G$ on $\K_{\M,f}^{\ssb}$ and
$\K_{\M,f}^{\ssb}[\dt^{-1}]$ in the same way as (1.5.2).
\sk
To show the assertion (4.1.2), consider the {\it relative logarithmic}
de Rham complex
$$\A_{\M,f}^{\prime\,\ssb}:=\M\sotim_{\OX}\Omega_{X/S}^{\ssb}(\log E),$$
where
$$\Omega_{X/S}^i(\log E):=\Omega_X^i(\log E)/
\bl(df/f\wedge\Omega_X^{i-1}(\log E)\br).$$
There is an isomorphism of complexes
$$df/f\wedge:\A_{\M,f}^{\prime\,\ssb}\simto\A_{\M,f}^{\ssb}[1],
\leqno(4.1.3)$$
by using the acyclicity of the complex
$\bl(\Omega_X^{\ssb}(\log E),df/f\wedge\br)$ together with the relation
$$d(df/f\wedge\omega)=-df/f\wedge d\omega\q\h{for}\,\,\,
\omega\in\Omega_X^i(\log E).$$
The following is known to the specialists:
\ms\nin
(4.1.4)\q The action of $t\dt$ on $\Hc^i\A_{\M,f}^{\prime\,\ssb}$
corresponds to that of $\dt t=t\dt+1$ on $\Hc^{i+1}\A_{\M,f}^{\ssb}$
under the isomorphism induced by (4.1.3).
\ms
Indeed, the action of $t\dt$ on $\Hc^i\A_{\M,f}^{\prime\,\ssb}$ and
that of $\dt t$ on $\Hc^{i+1}\A_{\M,f}^{\ssb}$ can be defined by using
the boundary morphism of the long exact sequence associated with the
short exact sequence of complexes
$$0\to\A_{\M,f}^{\ssb}\to\M\sotim_{\OX}\Omega_X^{\ssb}(\log E)
\to\A_{\M,f}^{\prime\,\ssb}\to 0,$$
together with the isomorphism (4.1.3). More explicitly, we have for
$[\omega']\in\Hc^i\A_{\M,f}^{\prime\,\ssb}$
$$t\dt[\omega']=[\eta']\,\,\,\h{if}\,\,\,
(df/f)\wedge\eta'=\nabla\omega',$$
and for $[\omega]\in\Hc^{i+1}\A_{\M,f}^{\ssb}$
$$\dt t[\omega]=[\nabla\eta]\,\,\,\h{if}\,\,\,
(df/f)\wedge\eta=\omega,$$
where $[\omega]$ denotes the class of $\omega$, etc. So we get (1.4.4).
(These actions are compatible with the $\D_{S,0}$-module structure
on the Gauss-Manin systems $\G_{\M,f}^i$ by the canonical morphism
$\A_{\M,f}^{\ssb}\into\K_{\M,f}^{\ssb}$.)
Note that $df/f\wedge$ (instead of $df\wedge$) is used in the above
argument, and this gives the difference between $t\dt$ and $\dt t$,
and also $[0,1)$ and $[-1,0)$ before (4.1.2).
\sk
The proof of (4.1.2) for $\Hc^i\A_{\M,f}^{\prime\,\ssb}$ (with $-1$
replaced by $0$) can be reduced to the case ${\rm rank}\,\M=1$ by
increasing induction on ${\rm rank}\,\M$ using a long exact
sequence associated with a short exact sequence of locally free
sheaves with regular singular integrable connection, since the
assertion is local. In case ${\rm rank}\,\M=1$, we can calculate
$\Hc^i\A_{\M,f}^{\prime\,\ssb}$ explicitly by using a Koszul
complex as in [St], see (4.2) below.

\ms\nin
{\bf 4.2.~Rank 1 case.} With the above notation and the assumption,
assume ${\rm rank}\,\M=1$. Let $x_1,\dots,x_n$ be local coordinates
of $X$ such that we have locally
$$f=\mprod_{i\in J}\,x_i^{e_i},\q E=\mcup_{i\in J'}\,\{x_i=0\},$$
where $J\subset J'$ and $e_i\ges 1$ for $i\in J$. Then the 
logarithmic de Rham complex $\M\sotim_{\OX}\Omega_X^{\ssb}(\log E)$
is locally the Koszul complex associated with the operators $\xi_i$
on $\OX$ for $i\in[1,n]$, defined by
$$\xi_i=\begin{cases}x_i\dd_{x_i}+\al_i&\h{if}\,\,\,i\in J',\\
\dd_{x_i}&\h{if}\,\,\,i\notin J',\end{cases}$$
where $\al_i\in\C$ and its real part is contained in $[0,1)$.
The relative logarithmic de Rham complex
$\A_{\M,f}^{\prime\,\ssb}=\M\sotim_{\OX}\Omega_{X/S}^{\ssb}(\log E)$
is the quotient complex defined by the relation
$$\msum_{i=1}^{r'}\,e_i\,[dx_i/x_i]=0,$$
where $[dx_i/x_i]$ denotes the class of $dx_i/x_i$. Choose $i_0\in J$.
Then $\A_{\M,f}^{\prime\,\ssb}$ is locally the Koszul complex associated
with the operators $\xit_i$ on $\OX$ for $i\in[1,n]\setminus\{i_0\}$,
defined by
$$\xit_i=\begin{cases}(x_i\dd_{x_i}+\al_i)-
(e_i/e_{i_0})(x_{i_0}\dd_{x_{i_0}}+\al_{i_0})&\h{if}\,\,\,
i\in J\setminus\{i_0\},\\
x_i\dd_{x_i}+\al_i&\h{if}\,\,\,i\in J'\setminus J,\\
\dd_{x_i}&\h{if}\,\,\,i\notin J'.\end{cases}$$
This implies that $\Hc^p\bl(\A_{\M,f}^{\prime\,\ssb}\br)_0=0$ in case
$\al_i\ne 0$ for some $i\in J'\setminus J$, since
${\rm Re}\,\al_i\in[0,1)$. We may thus assume
$$\al_i=0\,\,\,\h{for any}\,\,\,i\in J'\setminus J.$$
Set
$$A_0:=\bl\{\msum_{\nu\in\NN^J}\,a_{\nu}x^{\nu}\in\C\{x_J\}
\,\big|\,\nu_i+\al_i=(e_i/e_{i_0})(\nu_{i_0}+\al_1)\,\,
(\forall\,i\in J\setminus\{i_0\})\br\},$$
where $\C\{x_J\}:=\C\{x_i\,(i\in J)\}$ and
$\nu=(\nu_i)_{i\in J}\in\NN^J$.
We then get
$$\Hc^p\bl(\A_{\M,f}^{\prime\,\ssb}\br)_0=\mopl\,A_0\,
[dx_{i_1}/x_{i_1}]\wedge\cdots\wedge[dx_{i_p}/x_{i_p}],$$
where the direct sum is taken over
$\{i_1,\dots,i_p\}\subset J\setminus\{i_0\}$ with
$i_1<\cdots<i_p$. (This is a generalization of [St] which treated
the constant local system case.)
\sk
We can calculate the action of $t\dt$ on
$$\C\,x^{\nu}\,[dx_{i_1}/x_{i_1}]\wedge\cdots\wedge
[dx_{i_p}/x_{i_p}]\subset\Hc^p\bl(\A_{\M,f}^{\prime\,\ssb}\br)_0,$$
with $x^{\nu}\in A_0$. By the argument in (4.1),
this is given by the multiplication by
$$(\nu_i+\al_i)/e_i,$$
which is independent of $i\in J$ by the definition of $A_0$.
\sk
The localization of $\Hc^p\bl(\A_{\M,f}^{\prime\,\ssb}\br)_0$ by $t$
is the localization of $\A_{\M,f}^{\prime\,\ssb}$ by $f$,
and the above calculation can be extended to the localization.
Here the real part of $(\nu_i+\al_i)/e_i$ is strictly negative
in case $\nu_i\in\Z_{<0}$, since ${\rm Re}\,\al_i\in[0,1)$.
This implies that $\Hc^p\bl(\A_{\M,f}^{\prime\,\ssb}\br)_0$ is a finite
free $\C\{t\}$-module which is stable by the action of $t\dt$ and
such that the real part of any eigenvalue of the residue of $t\dt$ is
contained in $[0,1)$. This finishes the proof of (4.1.2).

\ms\nin
{\bf 4.3.~Proof of Theorem~3.} The canonical morphisms of Theorems~2
and 3 are constructed in the same way as the local system case. The
assertions are then reduced to the case $\M$ is a regular holonomic
$\D$-module. Let $E$ be a divisor on $X$ such that the restriction
of $\DR(\M)$ to the complement of $E$ is a local system supported
on a smooth variety. Let $i:E\into X$ and $j:X\setminus E\into X$
denote the inclusions. We have a distinguished triangle
$$i_*i^!\M\to\M\to\R j_*j^*\M\to.$$
By increasing induction on $\dim{\rm supp}\,\M$, the assertions
are then reduced to the case where $\M=\R j_*j^*\M$ (using the
canonical morphisms mentioned above).
Here we may assume that $\M$ is supported on $X$ and $D\cup E$ is a
divisor with simple normal crossings on $\X$ by replacing $\X$ with
a desingularization of the closure of the support of $\M$ in $\X$ if
necessary (using the compatibility with the direct image explained at
the end of (2.1)).
Note that the assertions are easy to show in case the support of
$\M$ is contained in a fiber of $f$.
Replacing $E$ with a larger divisor and blowing-up further $\X$ if
necessary, we may moreover assume:
\ms\nin
(4.3.1)\q $\fb({\rm Sing}\,\fb)\cap f(E)=\emptyset$.
\ms\nin
(4.3.2)\q If $E$ contains a divisor $Z$ of $X$ with $f(Z)$ a
point, then $E$ contains $f^{-1}f(Z)$.
\ms
Then Theorem~3 follows since the proofs of Propositions~2, 3, and 4
in the local system case are generalized to this situation by using
the arguments in (4.1-2).
(For instance, in the generalization of (2.4) with ${\rm rank}\,\M=1$,
there are three subsets $J,J',J''$ of $[1,n]$ respectively
corresponding to $f$, $D$, and the closure of $E$, where
$J'\cap J''=\emptyset$ and $J\subset J'\cup J''$ by (4.3.2).
If $J\ne\emptyset$ (i.e. if $f\ne 0$), then we can reduce to the case
$J=J'\cup J''=[1,n]$ with $J'$ non-empty by the same argument as in
(2.4). The argument for the case $J=\emptyset$ is similar and easier.
The details are left to the reader.)
This finishes the proof of Theorem~3.

\end{document}